\input amstex
\input amsppt.sty

\NoBlackBoxes

\def\epsilon{\varepsilon}

\def\phi{\varphi}

\document

\topmatter
\title $L_h^2$-domains of holomorphy in the class of unbounded Hartogs domains
\endtitle
\author Peter Pflug (Oldenburg) and W\l odzimierz Zwonek (Krak\'ow)
\endauthor
\abstract A characterization of $L_h^2$-domains of holomorphy in the class of Hartogs domains in $\Bbb C^2$ is given.
\endabstract
\address
Carl von Ossietzky Universit\"at Oldenburg, Institut f\"ur
Mathematik, Postfach 2503, D-26111 Oldenburg, Germany
\endaddress
\address
Uniwersytet Jagiello\'nski, Instytut Matematyki, Reymonta 4,
30-059 Krak\'ow, Poland
\endaddress
\email pflug\@mathematik.uni-oldenburg, Wlodzimierz.Zwonek\@im.uj.edu.pl
\endemail
\leftheadtext{Peter Pflug \& W\l odzimierz Zwonek}
\rightheadtext{$L_h^2$-domains of holomorphy}
\thanks Mathematics Subject Classification (2000). Primary: 32A07, 32D05.
%Secondary: 32U35, 30C85.
\endthanks
\thanks The research was supported by DFG Project No. 227/8-1 and was a part of the Research Grant
 No. 1 PO3A 005 28, which is
supported by public means in the programme promoting science in
Poland in the years 2005-2008. The second author is a fellow of
Rector's Scholarship Fund at the Jagiellonian University
\endthanks
\endtopmatter

There is a precise geometric characterization of bounded $L_h^2$-domains of holomorphy. Namely, we have the following

\proclaim{Theorem 1 {\rm (see \cite{Pfl-Zwo~1})}} Let $D$ be a bounded pseudoconvex domain in $\Bbb C^n$. Then $D$ is
an $L_h^2$-domain of holomorphy if and only if $U\setminus D$ is not pluripolar
for any open set $U$ with $U\setminus D\neq\emptyset$.
\endproclaim

As noted by M. A. S. Irgens there is no hope that an analogous result will hold in the unbounded case
-- it is sufficient to take the domain $\Bbb C\times\Bbb D$ (see \cite{Irg}). It is not an $L_h^2$-domain of holomorphy
(the space $L_h^2(\Bbb C\times\Bbb D)$ is trivial),
although the geometric condition from Theorem 1 is satisfied.

Therefore, it is natural to try to find a characterization of unbounded $L_h^2$-domains of holomorphy. Recall that
there is such a characterization in the case of planar domains.

\proclaim{Theorem 2 {\rm (see e.g. \cite{Con}, Chapter 21.9)}} Let $D$ be a domain in $\Bbb C$.

Then $D$ is an $L_h^2$-domain of
holomorphy if and only if $U\setminus D$ is not polar
for any open set $U$ with $U\setminus D\neq\emptyset$.

More precisely, for a point $a\in\partial D$ and an open neighborhood $U$ of $a$
there is an analytic continuation of any function $f\in L_h^2(D)$ onto $U$ if and only if $U\setminus D$ is polar.
\endproclaim

Another class of domains in which a full description of $L_h^2$-domains of holomorphy is known
is the class of Reinhardt domains (see \cite{Jar-Pfl}).

In the paper we present a characterization of $L_h^2$-domains of holomorphy in the class of unbounded Hartogs
domains with the base being a planar domain. The results of the paper may also be seen as a continuation
of results from \cite{Pfl-Zwo~2}, where Bergman completeness in the class of unbounded Hartogs domains is studied.

\bigskip

For a domain $D\subset\Bbb C^n$ denote by $L_h^2(D)$ the class of square integrable holomorphic functions on $D$.

Recall that a domain $D\subset\Bbb C^n$ is called an {\it $L_h^2$-domain of holomorphy } if there are no
domains $D_0,D_1\subset\Bbb C^n$ with $\emptyset\neq D_0\subset D_1\cap D$, $D_1\not\subset D$ such that
for any $f\in L_h^2(D)$ there exists an $\tilde f\in\Cal O(D_1)$ with $\tilde f=f$ on $D_0$.

\bigskip

For a subharmonic function $\rho:D\mapsto[-\infty,\infty)$, where $D$ is a domain in $\Bbb C$, define
$$
G_{D,\rho}:=\{(z_1,z_2)\in D\times\Bbb C:|z_2|<e^{-\rho(z_1)}\}.
$$
The above defined domain $G_{D,\rho}$ is a {\it Hartogs domain with the base equal to $D$}.

For $f\in L_h^2(G_{D,\rho})$ define
$$
\rho_f(z_1):=\limsup\sb{j\to\infty}\frac{1}{j}\log|f_j(z_1)|,\;
z_1\in D,
$$
where $f(z_1,z_2)=\sum\sb{j=0}\sp{\infty}f_j(z_1)z_2^j$,
$(z_1,z_2)\in D$; $f_j$'s are the coefficients of the Hartogs
expansion of $f$ in $G_{D,\rho}$. Certainly, $\rho_f\leq\rho$ on
$D$, so $\rho_f^*\leq\rho$ on $D$, where $g^*$ denotes the upper regularization of the function $g$.

Then define $\tilde\rho:=\sup\sb{f\in L_h^2(G_{D,\rho})}\rho_f^*$
on $D$. Certainly, $\tilde\rho^*$ is a subharmonic function on
$D$, $\tilde\rho^*\leq\rho$.

For a domain $D\subset\Bbb C$ define
$$
\Cal S:=\Cal S(D):=\{z\in\partial D:U\setminus D \text{ is polar
for some open neighborhood $U$ of $z$}\}.
$$

Now we may formulate Theorem 2 as follows. The domain $D\subset\Bbb C$ is an $L_h^2$-domain of holomorphy iff
$\Cal S=\emptyset$.

We also denote by $\Bbb D$ the unit disc in $\Bbb C$.

Main aim of our paper is the following

\proclaim{Theorem 3} Assume that $\rho$ is bounded from below.

{\rm (a)} If $D\neq\Bbb C$ then $G_{D,\rho}$ is an $L_h^2$-domain of
holomorphy iff $\limsup\sb{D\owns z\to z_0}\rho(z)=\infty$ for any
$z_0\in \Cal S$.

{\rm (b)} If $D=\Bbb C$ then $G_{\Bbb C,\rho}$ is an $L_h^2$-domain of
holomorphy iff $\rho$ is not constant.
\endproclaim

Let us begin with some lemmas.

\proclaim{Lemma 4} Assume that $\rho$ is bounded from below and let
$z_1^0\in D$ be such that there is no open connected neighborhood
$U\subset D$ of $z_1^0$ such that $\rho_{|U}$ is constant. Then
$\rho(z_1^0)=\tilde\rho^*(z_1^0)$.
\endproclaim
\demo{Proof of Lemma 4} Assume that $\rho\geq M$ on $D$. Suppose
that $\tilde\rho^*(z_1^0)<M_1<\rho(z_1^0)$. Therefore, there is an
open disc $U$ with the center at $z_1^0$, $U\subset D$ such that
$\rho_f(z_1)\leq\rho_f^*(z_1)\leq\tilde\rho(z_1)<M_1$ for any
$z_1\in U$ and $f\in L_h^2(G_{D,\rho})$. Therefore, for any $f\in
L_h^2(G_{D,\rho})$ the function $\Cal F_f$ defined as
$$
\Cal F_f(z_1,z_2):=\sum\sb{j=0}\sp{\infty}f_j(z_1)z_2^j,\;
(z_1,z_2)\in U\times e^{-M_1}\Bbb D,
$$
is a well defined holomorphic function -- it follows from the
Hartogs Lemma that the Hartogs series defining $\Cal F_f$ is
locally uniformly convergent in $U\times e^{-M_1}\Bbb D$.
Certainly, $\Cal F_f=f$ on $G_{D,\rho}\cap(U\times e^{-M_1}\Bbb
D)$. Therefore, $\Cal F_f$ is a holomorphic continuation of $f$
onto $U\times e^{-M_1}\Bbb D$.

On the other hand let $z_2^0\in\Bbb C$ be such that
$|z_2^0|=e^{-\rho(z_1^0)}$. Note that $(z_1^0,z_2^0)\in(\partial
G_{D,\rho}\cap(U\times e^{-M_1}\Bbb D))$. For $z_2\in e^{-M_1}\Bbb
D$ define
$$
U(z_2):=\{z_1\in U:(z_1,z_2)\in G_{D,\rho}\}.
$$
We claim that there is a $\tilde z_2\in e^{-M_1}\Bbb D$ such that
$U(\tilde z_2)\neq\emptyset$ and $U\setminus U(\tilde z_2)$ is not
polar. Actually, since $\rho_{|U}$ is not constant we easily get the existence
of a $\tilde z_2\in e^{-M_1}\Bbb D$
such that $U(\tilde z_2)\neq\emptyset$ and $U(\tilde z_2)\neq U$.
Note that $U(\tilde z_2)=\{z_1\in
U:\rho(z_1)<-\log|\tilde z_2|\}$. Suppose that $U\setminus U(\tilde z_2)$ is polar.
Then $\rho(z_1)\geq -\log|\tilde z_2|$ for any $z_1\in U(\tilde z_2)$. Since $U\setminus U(\tilde z_2)$
is polar, $\rho(z_1)\geq -\log|\tilde z_2|$ for any $z_1\in U$, so
$U(\tilde z_1)=\emptyset$ -- contradiction. Therefore, $U\setminus U(\tilde z_2)\neq U$ is not polar.
In particular, there is a
function $f\in L_h^2(V_1\times\{\tilde z_2\})$, where
$V_1:=\{z_1\in D:(z_1,\tilde z_2)\in G_{D,\rho}\}$, which does not
have a holomorphic continuation on $U\times\{\tilde z_2\}$ (see
Theorem 2).

It follows from a result of T. Ohsawa (see \cite{Ohs}), applied to
$\Psi(\cdot):=2g_{e^{-M}\Bbb D}(\tilde z_2,\cdot)$ ($g_D(p,\cdot)$ denotes the Green function
of the domain $D$ with the logarithmic pole at $p$), that there is a
function $F\in L_h^2(G_{D,\rho})$ such that $F_{|V_1\times\{\tilde
z_2\}}\equiv f$ (here and later making use of the result of T. Ohsawa
we utilize a formulation from the paper \cite{Chen-Kam-Ohs}, p. 706).
But it follows from the earlier property that $F$
extends to a holomorphic function on $U\times e^{-M_1}\Bbb D$.
Therefore, $f$ extends holomorphically onto $U\times\{\tilde
z_2\}$ -- contradiction. \qed
\enddemo

\subheading{Remark} Let us make a remark on the proof of the above lemma. We
provided the proof with the help of a new extension result of T. Ohsawa. The result in \cite{Ohs} applies
to the unbounded case (unlike the one in the standard version of the extension theorem
in \cite{Ohs-Tak}) but there are some limits. Namely, the possibility of the extension of an $L_h^2$-function
from the hyperplane depends on the existence of a suitable plurisubharmonic function $\Psi$ -- in our proof this is the
Green function of the projection of the domain $G_{D,\rho}$ onto the second variable. At this place it is important
that the projection is bounded, in other words that the function $\rho$ is bounded from below.

\proclaim{Lemma 5} Let $\rho$ be bounded from below and not
constant and let $z_1^0\in D$ be such that there is an open
connected neighborhood $U\subset D$ of $z_1^0$ such that
$\rho_{|U}$ is constant. Then $\tilde\rho^*\equiv\rho$ on $U$.
%If,
%additionally, $\Bbb C\setminus D$ is not polar or $\rho$ is not
%constant then $\rho\equiv\tilde\rho$ on $U$. If $\Bbb C\setminus
%D$ is polar and $\rho$ is constant then $\rho\equiv-\infty$ on
%$D$.
\endproclaim
\demo{Proof of Lemma 5} Let $\tilde V_1$ denote the set of all
$z_1\in D$ such that there is an open neighborhood $V$ of $z_1$
such that $\rho_{|V}$ is constant. Certainly $\tilde V_1$ is open
and $\tilde V_1\supset U$. Denote by $V_1$ the connected component
of $\tilde V_1$ such that $V_1\supset U$. Then $\rho_{|V_1}\equiv
C$ for some $C\in\Bbb R$.

%It is sufficient to
%show that $\tilde\rho\equiv C$ on $V_1$.

We claim that
$$
\text{ $\tilde\rho^*$ is constant on $V_1$}.\tag{$\ast$}
$$
Without loss of generality assume that
$\tilde\rho^*\not\equiv-\infty$. Then
$(\tilde\rho^*)^{-1}(-\infty)$ is polar. We observe that in order to
prove \thetag{$\ast$} it is sufficient to show the following
property:

\bigskip

\item{($\ast\ast$)} For any $z_1\in V_1$ such that
$\tilde\rho^*(z_1)>-\infty$ there is an open neighborhood
$U(z_1)\subset V_1$ of $z_1$ such that $\tilde\rho^*$ is constant
on $U(z_1)$.

\bigskip

First we prove the implication
($\thetag{$\ast\ast$}\Rightarrow\thetag{$\ast$}$). Assume that
\thetag{$\ast\ast$} holds. Then the set $\{z_1\in
V_1:\tilde\rho^*(z_1)>-\infty\}$ is open, so
$(\tilde\rho^*)^{-1}(-\infty)\cap V_1$ is closed in $V_1$ and
polar. Therefore, the set $\{z_1\in
V_1:\tilde\rho^*(z_1)>-\infty\}$ is connected. But $\tilde\rho^*$
is locally constant there, so it is constant on
$V_1\setminus(\tilde\rho^*)^{-1}(-\infty)$. The subharmonicity of
$\tilde\rho^*$ implies that $\tilde\rho^*$ is constant on $V_1$.

\bigskip

Now we show the property \thetag{$\ast\ast$}.

Suppose that there is a $\tilde z_1\in V_1$ such that
$(\tilde\rho)^*(\tilde z_1)>-\infty$ and $\tilde \rho^*$ is not
constant on any neighborhood of $\tilde z_1$. Without loss of generality we may assume that
$\tilde\rho^*(\tilde z_1)<\rho(\tilde z_1)$. Let
$-\infty<M<\tilde\rho^*(\tilde z_1)$ and $M\leq \rho$ on $D$. The
function $\psi:=\max\{M,\frac{\rho+\tilde\rho^*}{2}\}$ defined on
$D$ is subharmonic, bounded from below, $\psi\leq\rho$,
$\tilde\rho^*(\tilde z_1)<\psi(\tilde z_1)<\rho(\tilde z_1)$ and
$\psi$ is not constant on any neighborhood of $\tilde z_1$.

Let $\tilde\psi$ denote the function defined for $\psi$
analoguously to the way the function $\tilde \rho$ was defined for
$\rho$. Note that $\tilde\psi^*\leq\tilde\rho^*$ on $D$, so
$\tilde\psi^*(\tilde z_1)\leq\tilde\rho^*(\tilde z_1)<\psi(\tilde z_1)$.
However, it follows from Lemma 4 applied to $\psi$ that
$\psi(\tilde z_1)=\tilde\psi^*(\tilde z_1)$ -- contradiction.

\bigskip

Consequently, \thetag{$\ast$} is satisfied, so
$\tilde\rho^*\equiv\tilde C\in[-\infty,\infty)$ on $V_1$.

We want to show that $\tilde\rho^*\equiv\tilde C=C\equiv\rho$ on
$V_1$. If $D\setminus V_1$ is polar then $\rho$ is constant on $D$
-- contradiction. Therefore, $D\setminus V_1$ is not polar, so
$\partial V_1\cap D$ is not polar, either. Therefore, there is a
point $\tilde z_1\in\partial V_1\cap D$ such that $V_1$ is not
thin at $\tilde z_1$, so $\tilde\rho^*(\tilde z_1)=\tilde C$;
moreover, $\rho$ is not constant on any neighborhood of $\tilde
z_1$, so in view of Lemma 4, $\tilde\rho^*(\tilde z_1)=\rho(\tilde
z_1)=C$, so $C=\tilde C$. \qed
\enddemo

As a consequence of Lemmas 4 and 5 we get the following result.

\proclaim{Corollary 6} Let $\rho$ be bounded from below and not
constant. Then $\rho=\tilde\rho^*$ on $D$.
\endproclaim

Let us define $\hat D$ to be the set of points from $D$ and those
points $\hat z_1\in\Cal S$ such that $\limsup\sb{D\owns
z_1\to\hat z_1}\rho(z_1)<\infty$. Note that $\hat D$ is a
domain with $D\subset\hat D\subset D\cup\Cal S$. We may
also define the function
$$
\hat\rho(\tilde z_1):=\cases \rho(z_1),& \text{ if $z_1\in D$},\\
\limsup\sb{D\owns z_1\to\hat z_1}\rho(z_1),& \text{ if
$z_1\in\tilde D\setminus D$}\endcases.
$$
Note that $\hat\rho$ is subharmonic on $\hat D$.

%Note that in the case of Theorem 4 $\Cal S=\Bbb C\setminus
%D=\partial D$.
\demo{Proof of Theorem 3} Let $\rho\geq M$ on $D$.

If $\rho$ is constant, $\Cal S=\emptyset$ and $\Bbb
C\setminus D$ is not polar then the domain $D$ is an $L_h^2$-domain of
holomorphy (see Theorem 2) and $G_{D,\rho}=D\times e^{-M_1}\Bbb D$ for some
$M_1\in\Bbb R$ ($M_1\equiv\rho$). Consequently, $G_{D,\rho}$ is
an $L_h^2$-domain of holomorphy.

First we show the sufficiency of the condition.

Assume that $\limsup\sb{D\owns z\to z_0}\rho(z)=\infty$ for any
$z_0\in\Cal S$ and that $\rho$ is not constant (in the case $\Cal
S\neq\emptyset$ the latter condition follows directly from the
first one). Suppose that $G_{D,\rho}$ is not the $L_h^2$-domain of
holomorphy. Then there are discs $P_j$, $Q_j$, $j=1,2$ such that
$P_j\subset\subset Q_j$, $j=1,2$, $P:=P_1\times P_2\subset
G_{D,\rho}$, $\partial G_{D,\rho}\cap\partial P\neq\emptyset$
and for any $f\in L_h^2(G_{D,\rho})$ there is a
$g\in \Cal O(Q_1\times Q_2)$ such that $f=g$ on $P_1\times P_2$.

Let us consider three cases.

Case I. $\partial P\cap\partial G_{D,\rho}\subset\Cal S\times\Bbb
C$. Then our assumption implies that there is a point
$(z_1^0,z_2^0)\in\partial P\cap\partial G_{D,\rho}$ such that
$z_1^0\in\Cal S$, $z_2^0\neq 0$ and $z_2^0\in P_2$. Consider the
set
$$
U:=\{z_1\in Q_1\cap D:(z_1,z_2^0)\in G_{D,\rho}\}
$$
Note that $U\neq\emptyset$ (because $P_1\subset U$). Note also
that $Q_1\setminus U$ is not polar. In fact, the assumption on the
boundary behaviour of $\rho$ implies that there is a point $\tilde
z_1\in Q_1\cap D$ such that $(\tilde z_1,z_2^0)\not\in
G_{D,\rho}$, so $\rho(\tilde z_1)\geq -\log|z_2^0|$. The
existence of only a polar set of such points would, however, lead
to contradiction with the mean value property of subharmonic
functions. Therefore, there is a function $f\in L_h^2(\tilde U)$,
where $\tilde U:=\{z_1\in D: (z_1,z_2^0)\in G_{D,\rho}\}$, which
does not have a holomorphic continuation on $Q_1$ (see
\cite{Con}). There is a function $F\in L_h^2(G_{D,\rho})$ such
that $F(z_1,z_2^0)=f(z_1)$, $z_1\in\tilde U$ (apply \cite{Ohs}
with $\Psi(z_1,z_2):=2g_{e^{-M}\Bbb D}(z_2^0,z_2)$). But such a
function $F$ has a holomorphic continuation on $Q_1\times Q_2$,
from which we conclude the existence of a holomorphic continuation
of $f$ on $Q_1$ -- contradiction.

Case II. $\partial P\cap\partial G_{D,\rho}\cap((\partial
D\setminus\Cal S)\times \Bbb C)\neq\emptyset$. The proof in this
case is similar to that in Case I. There is a point
$(z_1^0,z_2^0)\in\partial P\cap\partial G_{D,\rho}$ such that
$z_1^0\in\partial D\setminus \Cal S$ and $z_2^0\in P_2$ (but we
have no guarantee that we may assume additionally that $z_2^0\neq
0$). Consider the set
$$
U:=\{z_1\in Q_1\cap D:(z_1,z_2^0)\in G_{D,\rho}\}
$$
Note that $U\neq\emptyset$ (because $P_1\subset U$). Note also
that $Q_1\setminus U$ is not polar. In fact, this follows directly
from the fact that $z_1^0\in\partial D\setminus \Cal S$,
definition of $\Cal S$ and the inclusion $Q_1\setminus D\subset
Q_1\setminus U$. Therefore, there is a function $f\in L_h^2(\tilde
U)$, where $\tilde U:=\{z_1\in D: (z_1,z_2^0)\in G_{D,\rho}\}$,
which does not have a holomorphic continuation on $Q_1$ (see
Theorem 2). But then the function $F$ defined by the formula $F(z_1,z_2):=f(z_1)$, $(z_1,z_2)\in G_{D,\rho}$
is from the class $L_h^2(G_{D,\rho})$. But such a
function $F$ has a holomorphic continuation on $Q_1\times Q_2$,
from which we conclude the existence of a holomorphic continuation
of $f$ on $Q_1$ -- contradiction.

Case III. $\partial P\cap\partial G_{D,\rho}\cap(D\times\Bbb
C)\neq\emptyset$.

\comment If there were only assuming the opposite ther taking
$\tilde z_1
 then for
any bidisc $P\subset G_{D,\rho}$ such that $\partial P\cap\partial
G_{D,\rho}\neq\emptyset$ the inequality $\partial P\cap\partial
G_{D,\rho}\cap(D\times\Bbb C)\neq\emptyset$ holds. Therefore, it
is sufficient to show that one cannot extend $L_h^2$-holomorphic
functions through points from $\partial G_{D,\rho}\cap
(D\times\Bbb C)$. Suppose the opposite. Then there are discs
$P_j$, $Q_j$, $j=1,2$ such that $P_j\subset\subset Q_j$, $j=1,2$,
$Q_1\subset\subset D$, $P_1\times P_2\subset G_{D,\rho}$,
$\partial G_{D,\rho}\cap\partial(P_1\times P_2)\neq\emptyset$ and
for any $f\in L_h^2(G_{D,\rho})$ there is an $F\in \Cal
O(Q_1\times Q_2)$ such that $f=F$ on $P_1\times P_2$.
\endcomment

Denote by $M_1:=\sup\{|z_2|:z_2\in P_2\}$,
$M_2:=\sup\{|z_2|:z_2\in Q_2\}$. Then $0<M_1<M_2$. From our
assumption we conclude the existence of a point $\tilde z_1\in
\bar P_1\cap D$ such that $\rho(\tilde z_1)\geq-\log M_1$. On the
other hand the extension property implies that for any $f\in
L_h^2(G_{D,\rho})$ and for any $z_1\in P_1$ the inequality
$e^{-\rho_f(z_1)}\geq M_2$ holds, so $\tilde\rho^*(z_1)\leq-\log
M_2$, $z_1\in P_1$, implying the inequality $\tilde\rho^*(\tilde
z_1)\leq-\log M_2<-\log M_1$, which contradicts the equality
$\rho(\tilde z_1)=\tilde\rho^*(\tilde z_1)$ following from
Corollary 6.

\bigskip

Now we prove the necessity of the condition.

Recall that $\hat D\setminus D$ is a polar set. Therefore,
$(\hat D\setminus D)\cap G_{\hat D,\hat\rho}$ is pluripolar.
Since $L_h^2$-holomorphic functions extend through pluripolar sets
it is easy to see that $L_h^2(G_{D,\rho})=L_h^2(\ G_{\hat
D,\hat\rho})_{|G_{D,\rho}}$, which gives us the necessity of the
condition.

\qed
\enddemo

The following description of $L_h^2$-holomorphic hulls of domains
$G_{D,\rho}$ follows from Theorem 3 and its proof.

\proclaim{Corollary 7} Assume that $\rho$ is bounded from below.
Assume that $\Bbb C\setminus D$ is not polar or $\rho$ is not
constant. Then the $L_h^2$-holomorphic hull of $G_{D,\rho}$ equals
$G_{\hat D,\hat\rho}$. If $\Bbb C\setminus D$ is polar and
$\rho$ is constant then the $L_h^2$-holomorphic hull of
$G_{D,\rho}$ equals $\Bbb C^2$. \endproclaim

\subheading{Remark} The problem of a full understanding of the structure of $L_h^2$-domains of holomorphy
is far from being solved. For instance, a natural question whether we may remove the assumption of lower boundedness
of the function $\rho$ in Theorem 3 remains open. On the other hand the methods used in the paper may be easily
transferred to Hartogs domains with higher dimensional bases. However, because of the lack of the full description
of $L_h^2$-domains
of holomorphy in $\Bbb C^n$, $n\geq 2$, the results obtained in this case would be much more incomplete.
We think that to find
a complete characterization of $L_h^2$-domains of holomorphy in the class of Hartogs domains in two dimensional case
(as well as in higher dimensional case of Hartogs domains or even in the general case of unbounded domains),
a completely different approach to the problem should be found.

Some other problems remain also open. For instance, is it true that if $D$ is a pseudoconvex domain, $D$ is locally
an $L_h^2$-domain of holomorphy (which means that the geometric condition from Theorem 1 is satisfied) and
$L_h^2(D)\neq\{0\}$ then $D$ must be an $L_h^2$-domain of holomorphy. Another natural problem would be to find
a description of $L_h^p$-domains of holomorphy.

\bigskip

We finish the paper with presenting some sufficient condition for a pseudoconvex domain $D$
to have infinitely dimensional Bergman space $L_h^2(D)$. This gives a partial answer to the following problem.
Is there a pseudoconvex domain having finite dimensional but non-trivial Bergman space? A non-pseudoconvex example
of that type was  given in \cite{Wie}.

\proclaim{Proposition 8} Let $D\subset\Bbb C^n$, $D\neq\Bbb C^n$,
be an $L_h^2$-domain of holomorphy and let $\{\phi_j\}_{j\in J}$
be a complete orthonormal system in $L_h^2(D)$. Assume that there is an open set $U$ such that $U\cap D\neq\emptyset$,
$U\not\subset D$ and for any $j\in J$ the function $\phi_j$ has an analytic continuation onto $U$. Then
$\operatorname{dim}L_h^2(G_{D,\rho})=\infty$. In particular, any $L_h^2$-domain of holomorphy, which is balanced,
a Hartogs domain or a Laurent-Hartogs domain different from $\Bbb C^n$, has infinitely dimensional Bergman space.
\endproclaim
\demo{Proof} Suppose the contrary. Then $J$ is finite and $J\neq\emptyset$.
Since $D$ is an $L_h^2$-domain of holomorphy there is a function $f\in L_h^2(D)$, which does not have an analytic
continuation onto $U$. But $f=\sum\sb{j\in J}\lambda_j\phi_j$, where $\lambda_j\in\Bbb C$. Since $J$ is finite
$f$ has analytic continuation onto $U$ -- contradiction.
\qed
\enddemo

\enddocument